\documentclass[12pt,reqno]{article}
\usepackage{latexsym}
\usepackage{amssymb}
\usepackage{amsmath}
\usepackage[all]{xy}

\newtheorem{thm}{Theorem}[section]
\newtheorem{defin}{Definition}[section]
\newtheorem{prop}{Proposition}[section]
\newtheorem{lemma}{Lemma}[section]

\newcommand{\Xbar}{\bar {X}}

\newcommand{\Ric}{\text{Ric}}

\newcommand{\kahler}{K\"ahler     }
\newcommand{\KE}{K\"ahler-Einstein     }

\newcommand{\vol}{\text{Vol}}
\newcommand{\hk}{hyperk\"ahler     }

\newcommand{\beeq}{\begin{equation}}
\newcommand{\eeeq}{\end{equation}}
\newcommand{\beit}{\begin{itemize}}
\newcommand{\eeit}{\end{itemize}}
\newcommand{\bedes}{\begin{description}}
\newcommand{\eedes}{\end{description}}
\newcommand{\been}{\begin{enumerate}}
\newcommand{\eeen}{\end{enumerate}}

\def\ZZ {{\mathbb Z}}

\def\RR {{\mathbb R}}
\def\CC {{\mathbb C}}
\def\PP {{\mathbb P}}

\def\De{\Delta}
\def\Om{\Omega}
\def\Ga{\Gamma}

\def\si{\sigma}

\def\de{\delta}
\def\om{\omega}
\def\ga{\gamma}
\def\ve{\varepsilon}

\newenvironment{proof}{\medskip
\noindent{\bf Proof: }}{{\hfill$\square$}{\medskip}}

\title{Existence of complete \kahler Ricci-flat metrics on  crepant resolutions}

\author{Bianca Santoro\footnote{Partially supported by NSF grant DMS-1007155.}}
\date{}

\begin{document}
\maketitle

\begin{abstract}
In this note,
we obtain existence results for complete Ricci-flat \kahler metrics
on crepant resolutions of singularities of Calabi-Yau varieties.
Furthermore, for certain asymptotically flat Calabi-Yau varieties,
we show that the Ricci-flat metric on the resolved manifold
has the same asymptotic behavior as the initial variety.

\end{abstract}

\section{Introduction}
\label{s.introduction}

In 1978, Yau \cite{Y1} proved the Calabi Conjecture by showing the existence
and uniqueness of \kahler metrics with prescribed Ricci curvature on compact
\kahler manifolds.

The natural extension of this theorem to noncompact manifolds introduces new
subtleties.
 A first major advance in this direction was
the proof by Tian and Yau
 \cite{TY1} of a  non-compact version of
the Calabi Conjecture on quasi-projective manifolds that can be compactified
by the addition of a smooth, ample divisor.
Shortly after, in \cite{TY2}, the two authors extended their results to the case
where the divisor has multiplicity greater than $1$ and is allowed to have
orbifold-type singularities.

In \cite{Santoro},
the asymptotic behavior of the metrics constructed by Tian and Yau in \cite{TY1} was fully described.
 First, a sequence of
explicit \kahler metrics with special approximating properties is constructed.
Using those metrics as a starting point, we then work out the asymptotic expansion
of the metrics given in \cite{TY1}.

In a different direction towards the generalization of the Calabi-Yau problem to non-compact manifolds,
Joyce \cite{Joyce1} proved the existence and uniqueness
of Asymptotically Locally Euclidean Ricci-flat \kahler metrics in each \kahler class
of $\CC^n / \Gamma$, for $\Ga \subset SL(n, \CC)$ a subgroup acting freely on $\CC^n \setminus \{0\}$.
We point out that for the case $n=2$, Kronheimer \cite{Kron1} had already obtained this result via
a different approach (\hk quotients).

Recently, van Coevering \cite{VanC} proved that a crepant resolution of
a Ricci-flat \kahler cone \footnote{A Ricci-flat \kahler cone is a metric cone over a  Sasaki-Einstein manifold of positive scalar curvature.}
admits a Ricci-flat \kahler metric asymptotic to the cone metric in every compactly supported
\kahler class of the resolution.

The aim of this paper is to study the Calabi Conjecture on manifolds which are
crepant resolutions of complete \kahler spaces with isolated singularities, and
ultimately to construct further examples of Calabi-Yau manifolds with
prescribed behavior at infinity.

One of our main results is the following:

\begin{thm}
\label{t.existence}
Let $(M,\om)$ be an n-dimensional, complete Calabi-Yau variety, whose singularities are contained in a compact subset $K$ of $M$.

Suppose that $(M \setminus K, \om |_{M\setminus K})$ is a Ricci-flat manifold, with  bounded sectional curvature and
$(\alpha,\beta)$-polynomial growth.

Let $\pi: \widetilde{M} \rightarrow M$ be a crepant resolution of $M$.

Then, for $\ve$ sufficiently small, there exists a
complete Ricci-flat \kahler metric $\widetilde g$, with associated \kahler form $\widetilde \om$, on $\widetilde M$, where
$$
\widetilde \om \in \pi^* [\om] - \ve^2[E],
$$
where $[E] $ is the Poincar\'e Dual of the $(2n - 2)-$homology class
of the exceptional set $E$.

\end{thm}
The definition of $(\alpha,\beta)$-polynomial growth will be provided in Section \ref{s.background}.

Although the results in Theorem~\ref{t.existence} are interesting in themselves,
being a natural expansion of the works of Yau \cite{Y1}, Tian-Yau \cite{TY1} and
Joyce \cite{Joyce} on the Calabi Conjecture, one of their virtues is to provide new
 examples of complete Calabi-Yau manifolds.
This will be carefully explained
 in Section \ref{s.applications}.

The above mentioned  applications have as inspiration the $4$-dimensional case of
{\em gravitational instantons}, connected \hk manifolds, with great relevance
in quantum gravity and low-energy supersymmetric solutions of String Theory.
Since the only compact examples of gravitational instantons are flat tori
and $K3$-surfaces, we should seek for complete, non-compact examples, where
compactness is replaced by a suitable condition of decay of the
\hk metric to a flat metric at infinity, in the following way.

Let $(\Xbar, \partial \Xbar)$ be a connected, orientable compact $4$--manifold with
smooth boundary $\partial \Xbar$.
Then, $\Xbar$ can be decomposed as
$$
\Xbar = K \cup N,
$$
where $K$ is a compact set and $N$, the {\em neck}, is diffeomorphic to $\partial \Xbar \times \RR^+$.

We will assume that $\partial \Xbar$ is a smooth fibration over a base $B$ (also smooth), with fiber $F$, which
will be assumed to be a compact flat manifold.

In \cite{CK1} and \cite{Etesi}, the terminology used for the various cases
split according to the dimension of the fiber $F$: if 
$F$ has dimension zero and the \hk metric decays to
the Euclidean metric, we refer to the manifold as
ALE (asymptotically locally Euclidean); the 
ALF (asymptotically locally flat) manifolds are such the metric decays to a
circle bundle metric, that is, $\dim F  = 1$; finally, the case when the 
fiber is 2 (resp. 3) dimensional are referred to by ALG (resp. ALH).

In \cite{Kron1} and \cite{Kron2}, Kronheimer classified all $2-$dimensional
ALE \hk manifolds  as resolutions of $\CC^2 / \Gamma$, where
$\Gamma$ is either $A_k, D_k$ or $E_j$, $j = 6,7, 8$.
For the ALF case, $A_k$ and $D_k$ families of ALF metrics
have been constructed by Cherkis and Kapustin \cite{CK2}
and explicitly by Cherkis-Hitchin \cite{CH}. Using String  Theory motivations,
they conjecture that those
may be the only ALF examples of gravitational instantons.

For the ALG case, it has also been conjectured that the only
examples are $D_k$, $k = 0, \cdots, 5$, and $E_j$, $j = 6,7, 8$.

In a recent paper, Hein \cite{Hijo} produced new families of complete Ricci-flat ALG and ALH spaces by
removing anticanonical curves from rational elliptic surfaces. He conjectures that such families should
be exhaustive up to certain deformations which are described in the paper.

Our next result includes both $4$-dimensional and  higher dimensional analogues of gravitational instantons, but without the \hk constraint.
The following theorem
provides  more specific results about the decay of the  metrics in Theorem \ref{t.existence}
for the special
case of manifolds which we call {\em asymptotically locally flat of order $k$}. A
precise definition is provided in Section \ref{s.background}.

\begin{thm}
\label{t.ALFk}
Let $(M, \om)$ satisfy the hypotheses in Theorem \ref{t.existence}, and suppose in addition that
$M$ is
asymptotically locally flat of order $k$.

Then, the complete Ricci-flat metric  $\widetilde \om$ given by Theorem \ref{t.existence} satisfies

\begin{equation}
\label{e.thm}
\widetilde \om  =  \om + C \frac{\sqrt{-1}}{2\pi} \partial \bar \partial (\rho^{2+k - 2n}) +
\partial \bar \partial \Psi,
\end{equation}
where $\rho$ is a distance function \footnote{The notion of distance function will be defined in the next section.} on $\widetilde M$, and
$\Psi$ is a smooth function on  $\widetilde M$ such that
$\nabla^\ell \Psi = O(\rho^{\gamma - \ell})$, where $\gamma \in (1+k - 2n, 2+k-2n)$.
\end{thm}

This theorem says that the $ALF_k$ behavior at infinity is preserved under
our construction. We shall exploit this  fact in Section \ref{s.applications}.

\bigskip

\noindent
{\bf Acknowledgements:} The author would like to thank Prof. Gang Tian
for suggesting this topic of research, and to Prof. Mark Stern for
numerous useful suggestions.

\section{Background Material}
\label{s.background}

Throughout this paper, we say that $f$ is {\em of order $O(\phi)$} if there
exists a uniform constant $C$ such that $|f| \leq C |\phi|$ near infinity.

We should start by providing the definition of  Asymptotically Locally Flat manifolds.
Since the author was not able to  find
in the literature a precise definition of those objects other than in the
$4$-dimensional case, we set our own, inspired by the definition
of Asymptotically Locally Euclidean manifolds in \cite{Joyce} and the discussion in \cite{Etesi},
and the definition of $4$-dimensional ALF manifolds described by Minerbe in \cite{Minerbe}.

As in the $4$-dimensional case, our varieties, now of dimension $n$,
 will be of the form
$X = \Xbar \setminus \partial \Xbar$, where $\partial \Xbar$ will
be assumed to be a fibration over a smooth base manifold $B$, with fibers $F$.
We will keep the same notation as before, where $\Xbar$ is described as
$ K \cup N$, where $K$ is a (large) compact set that contains all the singularities 
of $\Xbar$, and $N$ is diffeomorphic to
$\partial \Xbar \times \RR^+$ via a map
\begin{equation}
\label{e.asympcoordinate}
\Pi: N \rightarrow \partial \Xbar \times \RR^+,
\end{equation}
to be referred to as an {\em asymptotic coordinate system}.
Let $r$ be a coordinate for the $\RR$-component of the neck $N$.

Let $p: \partial \Xbar \rightarrow B$ be the projection of $\partial \Xbar$ onto the base $B$ of the fibration, with
$k$-dimensional fiber $F$.
Let $g_B$ be a fixed metric on $B$, and let $g_F$ denote a flat metric on the fiber $F$.
We will abuse notation and also denote by $g_F$ the extension of the metric  on the fiber to the total
space $\partial \Xbar$ by means of a horizontal distribution.

\begin{defin}
Let $g$ be a  metric on $X$. We say that $(X, g)$ is  {\em Asymptotically Locally Flat 
of order $k$} (or $ALF_k$, for short), if the
metric $g$ is asymptotic to a metric of the form
$$
h  = dr^2 + r^{2} p^*(g_B) + g_F.
$$

Furthermore, in order to make the decay condition precise,
we require that the push-forward metric $\Pi_*(g)$ satisfies
\begin{equation}
\label{e.ALFkdecay}
\nabla^m(\Pi_*(g) - h) = O(r^{2+k-n-m})
\end{equation}
outside the compact portion $K$ of the manifold $X$.
\end{defin}

\begin{defin}
\label{d.distance}
If $(X,g)$ is $ALF_k$,  we say that $\rho: X \rightarrow [1,\infty)$
is a {\em distance function} on $X$ if given any asymptotic coordinate
system $\Pi: N\rightarrow \partial \Xbar \times \RR$, we have
$$
\nabla^m(\Pi_*(\rho) - r) = O(r^{1+k-n-m}).
$$
outside the compact set $K$.
\end{defin}

Note  that the condition on Definition \ref{d.distance} is independent
of the choice of asymptotic coordinate system.

Examples of $ALF_k$ manifolds include, for $k = 0$, the $ALE$ manifolds of Joyce \cite{Joyce}
(in particular, the $4$-dimensional example of $\CC^2$ with the Eguchi-Hanson metric),
and for $k = 1$ and dimension $4$,  the Multi-Taub-NUT metric of Hawking. One of the applications
of our work is to explain a way to obtain new examples of $ALF_k$ manifolds for any $k$.

\bigskip

We shall  now characterize the necessary regularity and volume growth for our results.
The following two definitions were introduced by Tian and Yau in \cite{TY1}.
\begin{defin}
Let $(M,g)$ be a Ricci-flat manifold, and let $\alpha$ and $\beta$ be positive constants.

We say that $M$ has {\em $(\alpha,\beta)$-polynomial growth}
if, for a fixed point $x_0 \in M$, the volume of the geodesic ball $B_1(x)$ of radius $1$ around a point $x \in M$
satisfies $\vol_g(B_1(x)) \geq C^{-1}(1 + \text{dist}_g(x,x_0))^{-\beta}$, and the geodesic
ball of radius $R$ satisfies $\vol_g(B_R(x_0)) \leq C R^\alpha$, where  $C = C(M)$ is a uniform constant.
\end{defin}

\begin{defin}
Let $(M,g)$ be a  \kahler manifold.
We say that $(M,g)$ is of {\em quasi-finite geometry of order
$\ell + \delta$} if there exist $r>0$,
$r_1> r_2 > 0$ such that for any $x \in M$, there exists a holomorphic
map
$$\Phi_x: U_x \rightarrow B_r(x)$$ such that
$B_{r_2} \subset U_x \subset B_{r_1}$, and
$ \Phi_x^*(g)$ is a \kahler metric on $U_x \subset \CC^n$
such that its metric tensor has derivatives up to order $\ell$
bounded and $\de$-Holder-continuously bounded.
\end{defin}

\noindent
{\bf Remark (Tian and Yau, \cite{TY1}):}
If a  \kahler manifold $(M, g)$ has its sectional curvature and
covariant derivative of the scalar curvature bounded, then
$(M, g)$ is of quasi-finite geometry of order $2 + 1/2$.

\section{Proof of Theorem \ref{t.existence}}
In this section, we shall explain the proof of our main theorem.

We remind the reader that
we are assuming that the complete variety $(M, g)$, which is a smooth 
manifold outside a compact set $K$, is such that 
$(M \setminus K,\om|_{M \setminus K})$ is a Ricci-flat \kahler manifold with bounded sectional curvature
and $(\alpha,\beta)$-polynomial growth.
Denote by  $\rho$  the distance function of the manifold $M$.
Let $\pi: \widetilde{M} \rightarrow M$ be a crepant resolution.

According to Hironaka \cite{Hironaka}, any resolution of singularities of a complex variety can
be obtained as a finite sequence of blow-ups of regular subvarieties. Therefore, in order to obtain a
positive-definite \kahler form on the resolution $\widetilde{M}$ of $M$ with the same asymptotics as the pull-back
$\pi^*(\om_g)$, it suffices to prove this claim for a single blow-up.
Also, there is no loss in generality in assuming that the blow-up is centered at a point $x \in M$ instead of a higher-dimensional
subvariety.

\begin{lemma}
Let $\text{Bl}_x: \widetilde{M}_x \rightarrow M$ be the blow-up of $M$ at the point $x \in M$.
Then, there exists a complete \kahler metric
$\om$ on $\widetilde M$ such that
$$
\om  = \text{Bl}_x^*(\om_g)
$$
outside a compact set $K \subset  \widetilde{M}_x$ which contains $x$.
\end{lemma}

\begin{proof}
The $(1,1)$-form $\text{Bl}_x^*(\om_g)$ is positive semi-definite, and its degeneracy locus is given by the
exceptional divisor $E = (\text{Bl}_x)^{-1}(x)$.

We want to define a Hermitian metric on the line bundle $[E]$ associated to $E$, and use its curvature form
to define a $(1,1)$-form on $\widetilde{M}_x$ which is positive-definite when restricted to $E$.

Let $\cal U \subset M$ be a small neighborhood of $x \in M$, and let $z$ be a local coordinate system centered at $x$.
The blow-up map $\text{Bl}_x$ is an isomorphism outside $\cal U$, and on $\cal U$ can be described in local coordinates as
$$
\widetilde{\cal U} = \text{Bl}_x^{-1}({\cal U}) = \{ ( z, \ell) \in {\cal U} \times \PP^{n-1}: z_j \ell_i = z_i \ell_j  \hspace{0.3cm} \forall i,j  \}.
$$
Using these coordinates, we define a metric on $[E]$ as follows: we use the above coordinates
to define a metric $h_1$ on $[E]|_{\widetilde {\cal U}}$ by
$$
h_1(z, \ell) = |(\ell_1, \cdots, \ell)| = ||\ell||^2
$$
Let $\si$ be a defining section of the exceptional divisor $E$. On $\widetilde M \setminus E$, we define the metric
$h_2$ by setting $\si(p) = 1$ for $p \in \widetilde M$.

Let $(\rho_1, \rho_2)$ be a  partition of unity subordinate to the cover $\{ \widetilde{{\cal U}_{2\ve}}, \widetilde M \setminus \widetilde{{\cal U}_{\ve}} \}$,
where ${\cal U}_\ve  = \{ |z| < \ve \} \subset {\cal U}$. The metric $h$ on $[E]$
is given by
$$
h = \rho_1 h_1 + \rho_2 h_2.
$$

Finally, our positive-definite \kahler form is given by
$$
\om = \text{Bl}_x^*(\om_g) + \alpha \frac{\sqrt{-1}}{2 \pi} \partial \bar \partial \log h,
$$
where $\alpha = \alpha(\ve, \rho_1) > 0$ is chosen to be small enough to ensure that the resulting
metric does not develop any degeneration.
Observe that outside $\widetilde{{\cal U}}_{2\ve}$, the second summand is identically zero, and
inside $\widetilde{{\cal U}_{\ve}}$, it is simply the pull-back of the Fubini-Study metric on $\PP^{n-1}$.
\end{proof}

The metric $\om$ on $\widetilde M$ given by the lemma satisfies
$
\Ric(\om) = 0
$
outside the compact set $K$, being identical to the pull-back of the original Ricci-flat metric
$\om_g$ on $M$.

Therefore, in order to establish the existence of a Ricci-flat \kahler metric on
$\widetilde M$, we need the following result.

\begin{thm}
\label{t.ALFMA}
Let $(X, g)$ be a  Calabi-Yau manifold with bounded geometry of order $2+ 1/2$ and of
$(\alpha,\beta)$-polynomial growth.
Let $f$ be a compactly supported function satisfying the integrability condition
$\int_X (e^f -1)\, \om^n = 0$.
Then, there exists a solution $u$ to the equation
\begin{equation}
\label{e.MA}
\begin{cases}
\left(\om + \frac{\sqrt{-1}}{2 \pi} \partial \bar \partial u\right)^n  =  e^f \, \om^n, \\
 \om + \frac{\sqrt{-1}}{2 \pi} \partial \bar \partial u  >  0  \hspace{1cm} \text {on $X$.}
\end{cases}
\end{equation}
\end{thm}

The proof of this theorem will follow closely the ideas on \cite{TY1}, with
some refinements provided by \cite{Santoro}.
The idea of the proof consists of obtaining the function $u$ as the uniform
limit of solutions $u_\ve$, as $\ve$ goes to zero,  of  perturbed Monge-Amp\`ere equations
\begin{equation}
\label{e.MAve}
\begin{cases}
\left(\om + \frac{\sqrt{-1}}{2 \pi} \partial \bar \partial u_\ve\right)^n  =  e^{f+ \ve u_\ve} \, \om^n, \\
 \om + \frac{\sqrt{-1}}{2 \pi}\partial \bar \partial  u_\ve  >  0  \hspace{1cm} \text {on $X$.}
\end{cases}
\end{equation}
The existence of those $u_\ve$, for every $\ve >0$, is
guaranteed by the main result (Theorem A) in \cite{CY1}, since we are assuming that our manifold
$(X,g)$ has bounded geometry.

The remainder of this section will be dedicated to the proof of uniform $C^{2,1/2}$-estimates
for $u_\ve$, since standard elliptic theory will guarantee that we can take the limit
when $\ve$  approaches zero and obtain a solution $u$ to (\ref{e.MA}).

We begin by checking that the integrability condition is satisfied for each $u_\ve$.
\begin{lemma}
\label{l.intcond}
Let $(X, g)$ be as in Theorem \ref{t.ALFMA}, and let $u_\ve$ be a solution to
(\ref{e.MAve}). Then
$$
\int\limits_M (e^{f + \ve u_\ve} - 1)  \, \om^n = 0
$$
\end{lemma}

\begin{proof}
In order to integrate by parts, we consider a cut-off function
$\Phi : \RR \rightarrow [0,1]$ such that $\Phi(t) = 1$ for $t \leq 1$,
and $\Phi(t) = 0$ for $t \geq 2$, and $|\Phi '| \leq 2$, and $|\Phi ''| \leq 4$.

Let $\rho(x)$ be the distance, with respect to the metric $g$, from $x$ to a fixed point
$x_0 \in X$. Set $\Phi_R(x) = \Phi({\rho(x)}/{R})$.

Multiplying both sides of (\ref{e.MAve}) by $\Phi_R(x)$, integrating by parts,
and using the fact
that all metrics $ \om_\ve := \om + \frac{\sqrt{-1}}{2 \pi} \partial \bar \partial u_\ve $
are equivalent to $\om$, it follows that
$$
\int\limits_M \left[\Phi_R(x)(e^{f + \ve u_\ve} - 1) \right]\, \om^n \leq \frac{C}{R} \int\limits_M |\nabla u_\ve|\, \om^n,
$$
where the derivative is taken with respect to the metric $g$.
Hence, the result amounts to showing that $\int\limits_M |\nabla u_\ve|\, \om^n$ is uniformly bounded.
Holder's inequality
gives us
$$
\int\limits_M |\nabla u_\ve|\, \om^n \leq \left(\int\limits_M \left[(1 + \rho(x))^{2q} |\nabla u_\ve|^2 \right] \, \om^n  \right)^{1/2}
 \left( \int\limits_M(1 + \rho(x))^{-2q} \, \om^n  \right)^{1/2},
$$
where $q$ is chosen to be $\alpha + 1$, where $\alpha$ is the volume growth rate of balls on the manifold $M$,
described in Theorem \ref{t.existence}.

Therefore, it only remains to prove that, for such $q$,
\begin{equation}
\int\limits_M \left[(1 + \rho(x))^{2q} |\nabla u_\ve|^2 \right]\, \om^n < \infty.
\end{equation}

This is a consequence of the following result:

\begin{lemma}
{\bf (\cite{TY1}, Lemma 3.3)}
For any constants $\ve >0$, $\ell \geq 1$ and $p\geq 0$,
$$
\int\limits_M \left([(1 + \rho)^p u_\ve]^\ell + |\nabla((1 + \rho)^p u_\ve) |^2  |(1 + \rho)^p u_\ve |^{2\ell - 2} \right)\, \om^n < \infty.
$$

\end{lemma}

For completeness, we provide a proof of this lemma in our simpler case ($\ell = 1$).
The idea is to use a cut-off function to obtain a recursive inequality
comparing the $L^2$-norm of $(1 + \rho)^p u_\ve$ and its derivative
with the $L^2$-norm of $(1 + \rho)^{p-1} u_\ve$. In finitely many  steps,
we can reach the critical value $-(\alpha + 1)$, and that will complete the proof.

We show that the inequality holds for $(u_\ve)_+ = \max \{ u_\ve, 0\}$ and $(u_\ve)_- =  -\min\{u_\ve, 0  \}$.

Let  $u = (u_\ve)_+$. Multiplying the Monge-Amp\`ere equation (\ref{e.MAve}) by
$(1 + \rho(x))^{2p}  \Phi_R^2(x) u $ (for simplicity, we will write $1 + \rho(x) = d$),
we obtain
\beeq
\label{e.lemma3.3}
 \int\limits_M \left[\left(d^{2p}  \Phi_R^2(x) u \right)(e^{f + \ve u} -1) \right]\, \om^n
=
\frac{\sqrt{-1}}{2\pi}  \int\limits_M \left[d^{2p}  \Phi_R^2(x) u \right]
\partial \bar \partial (u_\ve) \wedge \Om_\ve,
\eeeq
where $\Om_\ve = \left(\om^{n-1} + \cdots +  (\om + \frac{\sqrt{-1}}{2\pi} \partial \bar \partial (u_\ve))^{n-1} \right)$.

Rewriting $e^{f + \ve u_\ve} - 1$ as $e^f(e^{\ve u_\ve} -1) + e^f -1$,
integrating (\ref{e.lemma3.3}) by parts and regrouping the terms,
we obtain
\begin{multline}
\label{e.above1}
\int\limits_M \left[d^{2p}  (u \Phi_R^2(x))  e^f (e^{\ve u} -1) \right] \, \om^n
+
\int\limits_M \vert \nabla \left(d^{p} u \right) \vert^2  \Phi_R^2(x)  \, \om^n
\leq   \\ \leq
C \left\{
\int\limits_M \left[d^{2p}  u \Phi_R^2(x)  |e^{f} -1| \right]\, \om^n
+
\int\limits_M \left[d^{2p-2} u^2 \Phi_R^2(x)  \right]  \, \om^n +
% \right. + \\ + \left.
\int\limits_M \left[\frac{|\Phi'|^2(\rho/R)}{R^2} d^{2p} u^2 \right] \, \om^n
\right\},
\end{multline}
where the last term requires one further step of integration by parts.

Notice that the the first term on the left-hand side of (\ref{e.above1})
is bounded below by
$$
C \int\limits_M \left[ d^{2p}  \Phi_R^2(x) u^2  e^{f} \right] \, \om^n
\geq
C \int\limits_M \left[d^{2p}  \Phi_R^2(x) u^2 \right]   \, \om^n,
$$
since $e^{\ve u} - 1 \geq \ve u$, and $f$ has compact support.

The first integral on the right-hand side of (\ref{e.above1})
is finite, and $\frac{|\Phi'|(\rho/R)}{R^2} < C \Phi^2_{2R}$.
Hence, we have
\begin{equation}
\label{e.pinductive}
\int\limits_M \left[d^{2p}  \Phi_R^2(x) u^2 \right]   \, \om^n
+
\int\limits_M \left[ \vert \nabla(d^{p} u) \vert  \Phi_R^2(x) \right]  \, \om^n
\leq
C + C \int\limits_M \left [d^{2p-2}  \Phi_{2R}^2(x) u^2 \right] \, \om^n .
\end{equation}

Now, Theorem A in \cite{CY1} ensures that, for each $\ve >0$,
the solution $u_\ve$ is bounded. Our volume growth assumption ($M$
is of $(\alpha, \beta)$-polynomial growth) implies that
$$
\int\limits_M u^2 d^{-(\alpha + 1)}      \, \om^n  < \infty.
$$
Therefore, by inductively using (\ref{e.pinductive})  in order to
get $p$ to be negative enough, and of course letting $R \rightarrow \infty$,
we obtain the desired statement.

This completes the proof of the proposition.
\end{proof}

The higher order estimates are proved in \cite{Y1}:
\begin{thm}
(Yau, \cite{Y1})
There are constants $C_3$ and $C_4$, independent of $\ve$, such that
$$
0 \leq n + \Delta u_\ve \leq C_3 e^{C_4(u_\ve - \inf_X u_\ve)}.
$$

Furthermore, there are {\em a priori} bounds for \small{$|\nabla^3 u_\ve|$}
in terms of
\small{$\sup_X {|u_\ve|, |\Delta u_\ve|, |f|, |\nabla f|, |\nabla^2 f|} $}.
\end{thm}

Finally, the last ingredient on the proof of uniform bounds for $u_\ve$:
\begin{prop}
\label{p.sup}
There exists a constant $C$ (independent of $\ve$) such that
$\sup_X |u_\ve| \leq C$.
\end{prop}

The key observation about this proposition is that we can still use the weighted Sobolev inequalities
derived in \cite{TY1}, since our original manifold $(M, g)$ is Ricci-flat,
with the right regularity and volume growth.

\begin{proof}
The reader should be warned about the constant $C$ that will denote, unless
otherwise stated, various different constants which are all independent of $\ve$.

In the course of the proof, we shall make use of uniform bound on the Green's Kernel on annuli.

Applying the maximum principle for (\ref{e.MAve}), we can see that the maximum
and the minimum of $u_\ve$ are attained inside the support
of $f$. Let $x_{\text{max}}$ (resp. $x_{\text{min}}$) be points where the
maximum (resp. minimum) of $u_\ve$ are achieved.

The expression
$$\int\limits_M e^f(e^{\ve u_\ve}-1) \, \om^n = \int\limits_M \left(e^{f + \ve u_\ve}-1\right) \, \om^n + \int\limits_M \left(e^{f}-1\right) \, \om^n = 0
$$
implies that $u_\ve(x_{\text{max}}) > 0$ and  $u_\ve(x_{\text{min}}) < 0$.

Let us assume that the weighted average
$$
\text{Ave}_\rho(u_\ve) = \frac{\int\limits_M (1 + \rho(x))^{-N} u_\ve(x) \, \om^n}{\int\limits_M (1 + \rho(x))^{-N} \, \om^n}
$$
of $u_\ve$ with respect to the weight $(1 + \rho(x))^{-N}$ ($N = \alpha + 1$ chosen accordingly to the
volume growth of $\om$) is non-negative. The proof in the other case is analogous.

To complete the proof of this proposition, the following result is needed.

\begin{lemma}
\label{l.ave}
If $\text{Ave}_\rho(u_\ve) \geq 0$, then there exists a constant $C$ independent
of $\ve$ such that
$$
u_\ve (x_{\text{min}}) \geq -C \hspace{1cm}  \text{and} \hspace{1cm} \sup_M(u_\ve(x) - \text{Ave}_\rho(u_\ve) ) \leq C.
$$
\end{lemma}
Let us complete the proof of Proposition \ref{p.sup} assuming the Lemma \ref{l.ave},
and return to its proof later.

Consider a small convex geodesic ball $B_{2r}(x_{\text{min}})$ of radius $2r$ around $x_\text{min}$, and a cut-off
function $\eta$ supported on it, which is identically $1$ on $B_{r/2}(x_{\text{min}})$,
and zero outside $B_{3r/2}(x_{\text{min}})$.

Yau's higher order estimates \cite{Y1} give that
\begin{equation}
\label{e.yau}
\De u_\ve + n \leq C e^{C(u_\ve - \inf_M u_\ve)}.
\end{equation}
Set  $\psi = (u_\ve - \inf_M u_\ve - 1)_{-}$. Clearly, $0 \leq \psi \leq 1$,
$\psi(x_{\text{min}}) = 1$, and $\psi \neq 0$ imply that  $ u_\ve - \inf_M u_\ve < 1$.

Now, multiply (\ref{e.yau}) by $\eta^2(x) \psi(x) G(x_{\text{min}},x)$, where $G(x,y)$ is
the Green's kernel for the Dirichlet problem on the ball $B_r(x_{\text{min}})$, so as to obtain
$$
- \int\limits_M \De \eta^2(\psi)  \psi  G(x_{\text{min}},x) \, \om^n
\leq
C_1 \int\limits_M  \eta^2 \psi^2    G(x_{\text{min}},x) \, \om^n.
$$

Integrating the left-hand side by parts, we have
\begin{multline}
\label{e.blah}
\int\limits_M  \left[ \eta^2 |\nabla \psi|^2  G(x_{\text{min}},x)\right] \, \om^n
+
\int\limits_M \left[\nabla(\eta^2)(\psi \nabla \psi)  G(x_{\text{min}},x) \right]\, \om^n
+ \\ +
\int\limits_M \left[\eta^2 (\psi \nabla \psi) \nabla G(x_{\text{min}},x) \right]\, \om^n
\leq
C_1 \int\limits_M \left[  \eta^2 \psi^2  G(x_{\text{min}},x) \right] \, \om^n.
\end{multline}

Note that
$$
\int\limits_M \left[\eta^2 (\psi \nabla \psi) \nabla G(x_{\text{min}},x) \right] \, \om^n
= \frac{1}{2}\int\limits_M \left[\nabla(\eta^2 \psi^2) \nabla G(x_{\text{min}},x)\right] \, \om^n
- \int\limits_M \left[\nabla(\eta^2)\psi^2 \nabla G(x_{\text{min}},x) \right]\, \om^n,
$$
so we can rewrite the previous equation as
\begin{multline}
\label{e.blah2}
\int\limits_M  \left[ \eta^2 |\nabla \psi|^2  G(x_{\text{min}},x)\right] \, \om^n
+
\int\limits_M \left[\nabla(\eta^2)(\psi \nabla \psi)  G(x_{\text{min}},x) \right]\, \om^n
+
\frac{1}{2}\int\limits_M \left[\nabla(\eta^2 \psi^2) \nabla G \right]\, \om^n
\\ \leq
C_1 \int\limits_M \left[ \eta^2 \psi^2  G(x_{\text{min}},x) \right] \, \om^n +
\int\limits_M \left[\nabla(\eta^2)\psi^2 \nabla G(x_{\text{min}},x)\right] \, \om^n.
\end{multline}

Now, recalling that $G$ is non-positive on the support of $\eta$,
and using the inequality $2ab \leq a^2 + b^2$, we get that
$$\int\limits_M \nabla(\eta^2)(\psi \nabla \psi)  G(x_{\text{min}},x) \, \om^n \geq
\int\limits_M \psi^2 |\nabla \eta|^2 G(x_{\text{min}},x) \, \om^n + \int\limits_M  \eta^2 |\nabla \psi|^2  G(x_{\text{min}},x) \, \om^n.
$$

Therefore,
\begin{multline}
\frac{1}{2} \left(\int\limits_{B} \nabla(\eta^2 \psi^2) \nabla G(x_{\text{min}},x) \, \om^n
+
\int\limits_{B} |\nabla \psi|^2 \eta^2 G(x_{\text{min}},x) \, \om^n \right)
\leq \\
C_1
\int\limits_{B} \psi(x) G(x_{\text{min}},x) \, \om^n
+
\frac{1}{2}\left(\int\limits_{B} \psi^2 |\nabla \eta|^2 G(x_{\text{min}},x) \, \om^n
+
\int\limits\limits_{B} \psi^2 \nabla(\eta^2) \nabla G(x_{\text{min}},x) \, \om^n \right),
\end{multline}
where $B = B_{2r}(x_{\text{min}})$.

As  mentioned before, we want to use the fact that $G(x_{\text{min}},x)$
and $\nabla G(x_{\text{min}},x)$ are bounded independently of $\ve$ inside the annulus
$B_{3r/2}(x_{\text{min}}) \setminus B_{r/2}(x_{\text{min}})$, which contains the support of $\nabla(\eta)$.
We have, for $B=B_r(x_{\text{min}})$,
$$
1 = \psi^2(x_{\text{min}})  =
\int\limits_{B} \De  G(x_{\text{min}},x) \psi^2 \, \om^n \\
 \leq
C \left[ \int\limits_{B} \psi   G(x_{\text{min}},x)\, \om^n
+ \int\limits_{B} |\psi|^2 \, \om^n          \right],
$$
and Holder's inequality yields
\begin{eqnarray*}
1 \leq
C \left[\left(\int\limits_{B}|\psi|^{2n-2} \, \om^n \right)^{\frac{1}{2n-2}}
\left(\int\limits_{B}G(x_{\text{min}},x)^\frac{2n-1}{2n + 1} \, \om^n \right)^{\frac{2n-2}{2n-1}}
+ \int\limits_{B} |\psi|^2 \, \om^n\right].
\end{eqnarray*}

This implies that
there exists a constant $C>0$ independent of $\ve$ such that
$$C \, \text{vol}\left(\text{supp}(\psi) \cap {B_r(x_{\text{min}})}\right) \geq 1.$$
Hence, for $B = B_r(x_{\text{min}})$,
$$
\int\limits_M \left[(1 + \rho(x))^{-N} u_\ve \right] \, \om^n \leq
\int\limits_{M \setminus \text{supp}(\psi) \cap {B}} \left[(1 + \rho(x))^{-N} u_\ve \right] \, \om^n + C.
$$

However,
$$
\int\limits_{M \setminus \text{supp}(\psi) \cap {B}}  \left[(1 + \rho(x))^{-N} u_\ve \right]\, \om^n
\leq
\sup_M (u_\ve) \left(\int\limits_{M} (1 + \rho(x))^{-N} \, \om^n
\, - \, \int\limits_{\text{supp}(\psi) \cap {B}}(1 + \rho(x))^{-N} \, \om^n \right),
$$
and
$$\int_{ \text{supp}(\psi) \cap {B}}(1 + \rho(x))^{-N} \, \om^n \geq C(2 + \rho(x_{\text{min}}))^{-N}, $$
since the volume of the region of integration is bounded below.

Therefore,
$$
\int\limits_M (1 + \rho(x))^{-N} u_\ve \, \om^n \leq
C + \sup_M u_\ve \left[\int\limits_M (1 + \rho(x))^{-N} \, \om^n  - C(2 + \rho(x_{\text{min}}))^{-N}\right],
$$
which implies, together with the Lemma \ref{l.ave}, that
$\sup_M u_\ve \leq C$. This completes the proof of Theorem \ref{t.ALFMA} in the case where $\text{Ave}_\rho(u_\ve) \geq 0$.
The proof when $\text{Ave}_\rho(u_\ve) < 0$ is analogous, and will be omitted.
\end{proof}

\bigskip

\noindent
{\bf Proof of Lemma \ref{l.ave}:}

Let us write $\psi = (u_\ve - \text{Ave}_\rho(u_\ve))_+$, so as to have
$\psi (e^{\ve u_\ve} -1) \geq 0$ on $M$. We have
$$
- \int\limits_M |\nabla(\psi^{\frac{q+1}{2}})|^2 \, \om^n
=
- {\textstyle (\frac{q+1}{2})}^2 \int\limits_M \psi^{q -1} u_\ve (e^{f + \ve u_\ve}-1)\, \om^n.
$$

Thanks to (\ref{e.MAve}), the last estimate yields
$$
\int\limits_M  |\nabla(\psi^{\frac{q+1}{2}})|^2 \, \om^n \leq  {\textstyle (\frac{q+1}{2})}^2  \int\limits_M |\psi|^q |e^f -1| \, \om^n.
$$

Our manifold satisfies the $(\alpha,\beta)$-polynomial growth condition and is assumed to
have bounded sectional curvature, therefore
the weighted Sobolev inequalities developed in section $2$ of \cite{TY1} still hold. Therefore,
setting $d = 1 + \rho(x)$ as before, we have
\begin{equation}
\label{e.sobolevTY1}
\int\limits_M \left| \psi^{\frac{q+1}{2}} - \text{Ave}_\rho(\psi^{\frac{q+1}{2}})\right|^{\frac{2(2n+1)}{2n-1}} d^{-N} \, \om^n
\leq
C \int\limits_M |\psi|^q  d^{-N} \, \om^n,
\end{equation}
where the last constant depends on $q$, but not on $\ve$.

Since the weighted volume $\int\limits_M (1 + \rho(x))^{-N} \, \om^n = \int\limits_M d^{-N} \, \om^n $ is bounded
(for fixed  $N$ sufficiently large),
we can apply Holder's inequality again, to obtain (set $\ell = \frac{2n+ 1}{2n-1}$)
$$
\left\Vert (1 + |\psi|) \right\Vert_{\ell(q+1)} \leq C(q,n) \left\Vert (1 + |\psi|) \right\Vert_{q+1},
$$
where $\Vert f \Vert_p$ stands for the weighted $L^p-$norm
$ \Vert f \Vert_p = \left(\int\limits_M  f^p d^{-N}\, \om^n \right)^\frac{1}{p}.$

Set $q_0 = 2^\ell$, and $q_{j+1} = q_j^\ell$.
Inductively,
$$
\Vert (1 + |\psi|) \Vert_{q_{j+1}} \leq C(q,n) \Vert (1 + |\psi|) \Vert_{q_0} \leq
 C \int\limits_M  |\psi|\, d^{-N}\om^n
\leq C \left( \int\limits_M |\psi|^2  d^{-N} \, \om^n \right)^\frac{1}{2}
\left( \int\limits_M d^{-N} \, \om^n \right)^\frac{1}{2}.
$$

Note that the first term on the right-hand side of this inequality is bounded, as shown
in the proof of Lemma~\ref{l.intcond}. Therefore, by letting $j \rightarrow \infty$,
$$
\lim_{j \rightarrow \infty}|| \psi||_{q_j} = \sup_M(u_\ve - \text{Ave}_\rho(u_\ve))_+ \leq C,
$$
which completes the proof of the Lemma~\ref{l.ave}.

\vspace{1cm}

With the uniform bounds, we are able to consider a subsequence $u_\ve \rightarrow u$,
where $u$ is a solution to (\ref{e.MA}). We note that the solution $u$ will be bounded,
and $\int_X |\nabla u|^2 \, \om^n < \infty$.

To see this, we multiply  (\ref{e.MAve}) by $u_\ve$  and after integrating by parts
once, we conclude
$$
\int_X |\nabla u_\ve|^2 \, \om^n < C \int_X | u_\ve| |e^{\ve u_\ve } - 1|\, \om^n,
$$
which proves the claim.

\section{Asymptotics of the solution}

In this section, we shall give the proof of Theorem \ref{t.ALFk}.
The  main ideas are borrowed from
\cite{Santoro}, since we can take advantage of $f$ being compactly supported.
This theorem implies that the $ALF_k$-type of the manifold is
preserved under our construction.

We begin by showing the first order decay of $u$, which is a result that actually
is independent of  the topological type (at infinity) of our Ricci-flat manifold.
It only requires the same assumptions as for Theorem~\ref{t.existence}.

\begin{lemma}
\label{l.unifvanishing}
Let $(M, g)$ as in Theorem \ref{t.existence}, and let $u$ be the solution to (\ref{e.MA}) constructed
in Theorem \ref{t.existence}.
Then $u(x)$ converges uniformly to zero as $\rho(x) \rightarrow \infty$, where $\rho(x)$ is a distance function as
described in Definition \ref{d.distance}.
\end{lemma}

\begin{proof}
We shall use the fact that the solution $u$ was obtained as the uniform limit of solutions
$u_\ve$ to (\ref{e.MAve}), and prove uniform bounds on $u_\ve$ by using the Maximum principle for the complex
Monge-Amp\`ere operator $M: C^\infty(M,\RR) \rightarrow C^\infty(M, \RR)$, defined by
$$
M(\phi) = \log \left( \frac{(\om + \frac{\sqrt{-1}}{2 \pi}\partial \bar \partial \phi )^n}{\, \om^n}\right).
$$

Let $\rho$ be a distance function on the manifold $M$ (as in Definition \ref{d.distance}).
Observe that for $\beta \in (2, n-1)$, we have that
$M(C \rho^{-\beta}(x)) = C_\beta \rho^{-\beta -2}(x) + O(\rho^{-\beta -3}(x))$.
Also, for each solution $u_\ve$ of (\ref{e.MAve}), we have that $M(u_\ve) = f + \ve u_\ve$.

For a fixed $\delta >0$, set $C_1 = \frac{C_1'}{\de}$, where
$C_1' = \sup_{\{\rho(x) \geq \delta^{-1} \}} (|u_\ve| + 1)$, and set $C_2 = -C_1$.
We point out that, due to Proposition \ref{p.sup}, we can choose $C_1$ independently of $\ve$.

Hence, for all $x$ such that $\{ \rho^{-1}(x) = \de\}$, we have
$$
C_1 \rho^{-1}(x) > |u_\ve (x)| \hspace{1cm} \text{and} \hspace{1cm} C_2 \rho^{-1}(x) > |u_\ve (x)|.
$$
Furthermore, in the region $\{ x \in M;  \rho^{-1}(x) = \de\}$, for $\de$ chosen to be
sufficiently small, we have
$$
M(C_1 \rho^{-1}(x) ) < f + \ve u_\ve < M(C_2 \rho^{-1}(x) ).
$$

Finally, we observe that the solution $u_\ve$ (for each $\ve$) converges uniformly
to zero at infinity. Hence, we can apply the maximum principle to the operator
$M(.)$ so as to  conclude that there exists a constant $C$ such that
$$
-C \rho^{-1}(x) \leq u_\ve(x) \leq C \rho^{-1} (x).
$$

We complete the proof of this lemma by observing that thanks to  the uniform estimates in \cite{Santoro},
Proposition 5.2,
we can take the constant $C$ to be independent of $\ve$. Then, passing to the limit when
$\ve \rightarrow 0$, the claim follows.
\end{proof}

\begin{prop}
\label{p.vanishing}
Let $(M,g)$ as in Theorem~\ref{t.existence} and assume  that
$M$ is an $ALF_k$-manifold.
Let $u$ be the solution of $(\ref{e.MA})$ constructed in Theorem \ref{t.existence}.

Then $u$ has order of decay of
$O(\rho^{2+k - 2n}(x))$, where $n$ is the complex dimension of the manifold $M$ and $f \in C^\infty(M, \RR)$
has compact support.

\end{prop}

\begin{proof}
We shall use once more the maximum principle for the Monge-Amp\`ere operator $M(\phi)$, but
our choice of barrier shall be carried out in further detail.

The barrier function that will be used here is given by $\psi = C \rho^{2+k - 2n}(-\log \rho)^\ell$,
where the constants $C$ and $k$ will be carefully chosen in what follows.

After some computations, we obtain
$$ M (C \rho^{2+k - 2n}(-\log \rho)^\ell)  = C (2\ell + 1 + k - 2n)\rho^{k - 2n}(-\log \rho)^{\ell-1} \left( 1 + o(1)\right).
$$

Next, fix $\de > 0 $ .
and define
$C_1 = \frac{C_1'}{\de}$, where
$C_1' = \sup (|u| + 1)$, where the supremum is taken on the set $\{ x \in M;  \rho^{2+k - 2n}(-\log \rho)^\ell(x)  \leq \de  \}$.
Also, set $C_2 = -C_1$.
Clearly,
$$
C_1 \rho^{2+k - 2n}(-\log \rho)^\ell(x)  > u(x) \hspace{1cm} \text{and} \hspace{1cm} C_2 \rho^{2+k - 2n}(-\log \rho)^\ell(x)  < u(x).
$$ on the set where $\rho^{2+k - 2n}(-\log \rho)^\ell(x)  =\de$.

Now, we fix $\ell$ (the exponent of the log term) to be small enough so that $2\ell + 1+k  - 2n <0$.
Since $f$ has compact support \footnote{In fact, all we need from $f$ is that it decays strictly faster than $\rho^{-2n}$.},
there exists some small $\de> 0$ such that
\begin{eqnarray}
M(u) = f  & \geq & C_1 (2\ell+1+ k-2n) \rho^{-2n} (-\log \rho)^\ell(x)\left( 1 + o(1) \right) \\
          & = &  M (C_1 \rho^{2+ k  - 2n}(-\log \rho)^\ell(x)),
\end{eqnarray}
and analogously for the upper bound with $C_2$.

By Lemma \ref{l.unifvanishing}, we know that $u(x)$ converges uniformly to zero as $\rho(x) \rightarrow \infty$.

Therefore, we can apply the maximum principle to the operator $M$ and conclude that
there exists a constant $C$ such that
$$
|u (x)| \geq C \rho^{2 - 2n}(-\log \rho)^k(x).
$$
This completes  the proof of the proposition.
\end{proof}

\begin{prop}
\label{p.derivum}
Let $u$ be a solution to (\ref{e.MA}). Then, there exists
$C = C(k)$ such that, for all $x \in M$,
\begin{equation}
\label{e.derivum}
|\nabla^k u|_{g}(x) \leq C \rho^{\gamma - k} (x).
\end{equation}
\end{prop}

\begin{proof}
Since $f$ has compact support, the statement follows from observing that
the leading behavior of the derivatives of $u$ will depend solely on the term
$C \rho^{2+k - 2n}(x)$ appearing on the expansion of the solution $u$ given by Proposition \ref{p.vanishing}.
\end{proof}

\section{Application:  construction of complete non-flat Ricci-flat manifolds}
\label{s.applications}
\subsection{Asymptotically Locally Flat metrics on $\CC^n$, $n >2$}

In \cite{LeBrun}, LeBrun observed that $\CC^2$ admits a Ricci-flat \kahler metric
which is not flat. The so-called {\em Taub-NUT metric} of Hawking \cite{Hawking}
can be given explicitly on $S^3 \times \RR^+$ by
$$
g_T = \frac{\rho + 1}{4\rho} d\rho^2 + \rho (1 + \rho)[\sigma_1^2 + \sigma_2^2] + \frac{\rho}{\rho + 1}\sigma_3^2,
$$
where $\{ \sigma_1,\sigma_2,\sigma_3 \}$ is a left-invariant coframe for $S^3$ and $\rho \in \RR^+$.

Note that the volume of a large ball of geodesic radius $R$ about the origin is given by
$$
\text{Vol}_{g_{T}}(B_R(0)) = O(R^3),
$$
and this manifold is in fact an $ALF_1$-manifold \footnote{In the literature, an $ALF_1$ $4$-manifold is
called Asymptotically Locally Flat, or ALF for short}.
We also point out that the metric $g_t$ is in fact Hyperk\"ahler, since it has holonomy $SU(2)$.

For our example, we shall consider a copy of $\CC^2$, endowed with the Taub-NUT metric
$g_T$, and take the product with an Euclidean $(\CC^{n-2}, g_e)$. Clearly, the product
 metric $g \otimes g_e$ is a trivial example
of a non-flat complete Ricci-flat metric on $\CC^n$, but that is certainly not our main focus of interest.

The group $\ZZ^2$ acts in both components $\CC^2$ and $\CC^{n-2}$.
Let $\ga_1$ and $\ga_2$
be the generators of the action in each factor.
Consider the action on $\CC^n$ by the group (of order $2$)
$\Ga$ generated by the pair $(\ga_1, \ga_2)$.  The quotient $\CC^n/ \Ga$
is not a product of quotients, as $\Ga$ has order $2$.

Now, consider the manifold $M = (\CC^n, g_T \otimes g_e) / (\Gamma \times G)$. According to our construction,
$M$ will be a singular, complete Ricci-flat manifold, with a quotient singularity at the origin.

Let $\pi: \widetilde M \rightarrow M$ be a crepant resolution of $M$.
We shall point out here that the work of Sardo Infirri \cite{SI} implies the
existence of a crepant resolution for this quotient, since it is a toric variety.

Then, Theorem~\ref{t.existence}
can be applied to provide the existence of a complete, Ricci-flat  metric $\widetilde g$ on $\widetilde M$.
In view of Theorem~\ref{t.ALFk}, this metric will  also be a $ALF_1$ metric, since the
$ALF_k$-type is preserved under our construction.

\vskip 1cm

\flushleft

\medskip
{\bf Bianca Santoro} \ \  (bsantoro@ccny.cuny.edu)\\
The City College of New York, CUNY\\
Department of Mathematics \\
138th Street and Convent Avenue, NAC 6/203-C\\
New York, NY 10031\\
USA\\

\end{document}